\documentclass{amsart}
\usepackage[utf8]{inputenc}
\usepackage{amsmath}
\usepackage{mathtools}
\DeclarePairedDelimiter\floor{\lfloor}{\rfloor}

\theoremstyle{definition}
\newtheorem{theorem}{Theorem}[section]
\newtheorem{lemma}[theorem]{Lemma}
\newtheorem{definition}[theorem]{Definition}
\newtheorem{remark}[theorem]{Remark}
\newtheorem{proposition}[theorem]{Proposition}
\newtheorem{example}[theorem]{Example}

\title[Generalized Łojasiewicz's inequality \dots]{Generalized Łojasiewicz's inequality for definable and subanalytic multifunctions}
\author{Michał Kosiba}
\address{Jagiellonian University, Faculty of Mathematics and Computer Science, Łojasiewicza 6, 30-348 Kraków, Poland}
\email{michal.kosiba@student.uj.edu.pl}

\begin{document}

\begin{abstract}
This paper is devoted to proving the general Łojasiewicz inequality, in both the definable and subanalytic cases, under the most relaxed assumptions. It means that we drop the usual continuity and compactness assumptions. In the second part of the paper we concentrate on the Łojasiewicz inequality for multifunctions and apply it to the natural multifunctions related to the medial axis of a set (basic notion in pattern recognition).
\end{abstract}



\maketitle

\section{Introduction}
Among the strongest and most versatile tools in singularity theory, or analysis in the broad sense, one finds the so-called Łojasiewicz's inequalities. Generally speaking, they describe how two functions or sets approach one another and how this can be controlled under some (mild) conditions. The first type of Łojasiewicz's inequality was proved by prof. Łojasiewicz in the late '50s in order to solve Schwartz's problem of dividing a distribution by an analytic function. The aim of this paper is to formulate different types of Łojasiewicz's inequality both for functions and multifunctions. The results presented in this article can be divided into two parts.

The first part is devoted to the study of the Łojasiewicz inequality for functions. The most basic, classical version of this theorem (Theorem \ref{thm11} hereafter) is stated for two continuous, semialgebraic functions on compact subset. In \cite{pre} it is shown that under some conditions the assumption of continuity of one function can be omitted. In \cite{categ} we can find Łojasiewicz's inequality in any definable structure but in this case we require the continuity of both functions. Our goal will be to generalize these facts to any definable structure. We will show that under certain conditions we can omit the assumption of continuity not only of one functions but also of the other (Theorem \ref{uog} and Theorem \ref{glowne}). Moreover, we will prove that the functions can be defined on any definable set, not only compact or even closed. In the subanalytic case it turns out the situation is more subtle. Some kind of subanalytic versions (which work also in polynomially bounded definable structures) can be found in \cite{valette}. Although Theorem \ref{uog} has a direct subanalytic counterpart (Theorem \ref{subuog}), already with transposing Theorem \ref{glowne} we have to change the methods and add a boundedness assumption (Theorems \ref{subglowne} and \ref{subuog2}) which is unavoidable as shown in Ex 4.9. This is yet another case illustrating the rather big difference that exists between o-minimal structure and subanalytic geometry. 

The aim of the second part of the thesis is to use previously obtained results for definable multifunctions. At the beginning we will define several types of preimages of multifunctions that were presented in \cite{multi}. Then we will use them to define some counterparts of Łojasiewicz's inequality for definable multifunctions with compact values and later generalize it to the case of definable multifunctions with closed values. 
At the end of the thesis we will use the obtained inequalities for the multifunction of the closest points and the multifunction of the points realizing the distance to a set. These objects are connected with the notion of a medial axis which is a main concept of pattern recognition theory.

I would like to thank my supervisor Maciej Denkowski for suggesting the problem.

\section{Classic version of Łojasiewicz's inequality}
The facts and the definitions presented in this chapter come from \cite{coste}, \cite{long}, \cite{multi}. We gather here all the basic notions we will be using in the sequel.

\begin{definition}
A set $A \subset \mathbb{R}^{n}$ is called \textit{algebraic}, if there exists a subset $S \subset \mathbb{R}[X_{1},\ldots,X_{n}]$ such that
\[A=\{x \in \mathbb{R}^{n} \;|\;\forall\; f \in S: f(x)=0\}.\]
\end{definition}

\begin{definition}
A set $A \subset \mathbb{R}^{n}$ is called \textit{semialgebraic}, if it is of the form
\[A=\bigcup_{i=1}^{k} \bigcap_{l_{i}=1}^{j_{i}}\{P_{i}=0, Q_{il_{i}}>0\},\]
where $k,j_{i} \in \mathbb{N}$ and $P_{i},Q_{il_{i}} \in \mathbb{R}[X_{1},\ldots,X_{n}]$.
\end{definition}

\begin{definition}
Let $\mathcal{S}_{n}$ be a family of subsets of $\mathbb{R}^{n}$. We say that $\mathcal{S}_{n}$ is a \textit{Boolean algebra}, if:\\
1) $\emptyset \in \mathcal{S}_{n}$;\\
2) for $A,B \in \mathcal{S}_{n}$, we have $A \cap B$, $A \cup B$, $\mathbb{R}^{n} \backslash A$ $\in \mathcal{S}_{n}$.
\end{definition}

\begin{definition}
A \textit{structure} on the field $(\mathbb{R},+,\cdot)$ is a sequence $\mathcal{S}=\{\mathcal{S}_{n}\}_{n \in \mathbb{N}}$, such that for every $n \in \mathbb{N}$, $\mathcal{S}_{n}$ is a family of subsets of $\mathbb{R}^{n}$ satisfying:\\
1) $\mathcal{S}_{n}$ contains all the algebraic sets of $\mathbb{R}^{n}$;\\
2) $\mathcal{S}_{n}$ is a Boolean algebra of subsets of $\mathbb{R}^{n}$;\\
3) if $A \in \mathcal{S}_{n}$, $B \in \mathcal{S}_{m}$ then $A\times B \in \mathcal{S}_{n+m}$;\\
4) if $\pi: \mathbb{R}^{n}\times \mathbb{R} \longrightarrow \mathbb{R}^{n}$ is the natural projection and $A \in \mathcal{S}_{n+1}$, then $\pi(A) \in \mathcal{S}_{n}$.\\
The elements of $\mathcal{S}_{n}$ are called \textit{definable} sets.
\end{definition}

\begin{definition}
Let $\mathcal{S}$ be a structure on $(\mathbb{R},+,\cdot)$. We define a \textit{first-order formula} in the following way:\\
1) If $P\in \mathbb{R}[X_{1},\ldots,X_{n}]$, then $P(x_{1},\ldots,x_{n})=0$ and $P(x_{1},\ldots,x_{n})\geq 0$ are first-order formulas.\\
2) If $A \subset \mathbb{R}^{n}$ is a definable set, then $x \in A$ is a first-order formula.\\
3) If $\Phi(x_{1},\ldots,x_{n})$ and $\Psi(x_{1},\ldots,x_{n})$ are first-order formulas then
\[\Phi(x_{1},\ldots,x_{n})\vee \Psi(y_{1},\ldots,y_{m}),\; \Phi(x_{1},\ldots,x_{n})\wedge \Psi(y_{1},\ldots,y_{m}),\]
\[\Phi(x_{1},\ldots,x_{n}) \Longrightarrow \Psi(y_{1},\ldots,y_{m}), \;\sim \Phi(x_{1},\ldots,x_{n}),\; \exists \;x_{n}\in \mathbb{R}:\; \Phi(x_{1},\ldots,x_{n}),\]
\[\text{and} \;\forall\;x_{n} \in \mathbb{R}: \; \Phi(x_{1},\ldots,x_{n})\]
are all first-order formulas.
\end{definition}

\begin{theorem}[\cite{coste} Theorem 1.13]
Let $\mathcal{S}$ be a structure on  $(\mathbb{R},+,\cdot)$ and $\Phi(x_{1},\ldots,x_{n})$ a first-order formula $\mathcal{S}$. Then the set
\[\{(x_{1},\ldots,x_{n})\in \mathbb{R}^{n}: \; \Phi(x_{1},\ldots,x_{n})\}\]
is definable.
\end{theorem}

\begin{proposition}[\cite{coste} Proposition 1.12]
Let $A$ be a definable set. Then $int \; A$ and $\overline{A}$ are definable sets.
\end{proposition}

\begin{definition}
A structure $\mathcal{S}=\{\mathcal{S}_{n}\}_{n \in \mathbb{N}}$ on $(\mathbb{R},+,\cdot)$ is \textit{o-minimal}, if $\mathcal{S}_{1}$ consists of nothing else but all the finite unions of points and intervals of any type.

\end{definition}

\begin{theorem}[\cite{long} Theorem 1.6](Tarski-Seidenberg)
Let $\pi:\mathbb{R}^{n}\times \mathbb{R} \longrightarrow \mathbb{R}^{n}$ be the natural projection, $A \subset \mathbb{R}^{n+1}$ a semialgebraic set. Then $\pi(A)$ is semialgebraic.
\end{theorem}

\begin{proposition}[\cite{long} Proposition 1.4]
Any o-minimal structure contains semialgebraic sets.
\end{proposition}

Let us denote $\mathcal{SA}_{n}=\{A \subset \mathbb{R}^{n} \; | \; A \; \text{semialgebraic}\}$. 

\begin{theorem}[\cite{long}]
$\{\mathcal{SA}_{n}\}_{n \in \mathbb{N}}$ is an o-minimal structure on $(\mathbb{R},+,\cdot)$.
\end{theorem}

This theorem together with the previous proposition shows that any o-minimal structure contains semialgebraic sets as a substructure.

From this point we will understand definability as the definability in some o-minimal structure.

\begin{definition}
Let $A \subset \mathbb{R}^{n}$. The \textit{graph of a function} $f:A \longrightarrow \mathbb{R}^{m}$ is the set $\Gamma_{f}=\{(x,y)\in \mathbb{R}^{n}\times \mathbb{R}^{m} \;|\; f(x)=y\}$.
\end{definition}

\begin{definition}
Let $A \subset \mathbb{R}^{n}$. We say that a function $f:A \longrightarrow \mathbb{R}^{m}$ is \textit{definable} (in an o-minimal structure), if its graph is a definable set. 
\end{definition}

From the properties of the natural projection we get that the domain and the image of a definable function is a definable set.

\begin{proposition}
The restriction to a definable set and the composition of definable functions is again definable.
\end{proposition}

\begin{proposition}
\label{prop1}
Let $f,g:A \longrightarrow \mathbb{R}$, $A \subset \mathbb{R}^{n}$ be two definable functions.
Then the sets $\{x \in A \;|\; f(x)\leq g(x)\}$ and $\{x \in A \;|\; f(x) = g(x)\}$ are definable.
\begin{proof}[Proof] The functions $f$ and $g$ are definable, so the sets $\Gamma_{f}, \Gamma_{g}$ are definable too.
\[\{x \in A \;|\; f(x)\leq g(x)\}\ni z \Longleftrightarrow \]\[\Longleftrightarrow  \forall\; (x_{0},y_{0}) \in \Gamma_{f}, \forall\; (x_{1},y_{1})\in \Gamma_{g} :(z=x_{0} \wedge z=x_{1}) \Longrightarrow y_{0} \leq y_{1}.\]
\[\{x \in A \;|\; f(x)= g(x)\}\ni z \Longleftrightarrow \]\[\Longleftrightarrow \forall\; (x_{0},y_{0}) \in \Gamma_{f}, \forall \;(x_{1},y_{1})\in \Gamma_{g} :(z=x_{0} \wedge z=x_{1}) \Longrightarrow y_{0} = y_{1}.\]
We got the descriptions of both sets by first-order formulas which ends the proof.
\end{proof}
\end{proposition}

\begin{theorem}
\label{theo3}
Let $f,g:A \longrightarrow \mathbb{R}$, $A \subset \mathbb{R}^{n}$ be two definable functions.
Then $\{(x,y) \in \mathbb{R}^{n}\times\mathbb{R}^{n}\; |\; f(x) =g(y)\}$ and $\{(x,y) \in \mathbb{R}^{n}\times\mathbb{R}^{n}\; |\; f(x) \leq g(y)\}$ are definable sets.
\begin{proof}[Proof]
The sets $\Gamma_{f}, \Gamma_{g}$ are definable, so we get descriptions by first-order formulas:\\ 
\[\{(x,y) \in \mathbb{R}^{n}\times\mathbb{R}^{n}\; |\; f(x) =g(y)\}\ni (w,z) \Longleftrightarrow\] \[\Longleftrightarrow  \forall\; (x_{0},y_{0}) \in \Gamma_{f}, \forall\; (x_{1},y_{1}) \in \Gamma_{g}: (x_{0}=w \wedge x_{1}=z) \Longrightarrow y_{0}=y_{1}, \]
\[\{(x,y) \in \mathbb{R}^{n}\times\mathbb{R}^{n}\; |\; f(x) \leq g(y)\}\ni (w,z) \Longleftrightarrow \] \[\Longleftrightarrow  \forall\; (x_{0},y_{0}) \in \Gamma_{f}, \forall\; (x_{1},y_{1}) \in \Gamma_{g}: (x_{0}=w \wedge x_{1}=z) \Longrightarrow y_{0}\leq y_{1}. \]
\end{proof}
\end{theorem}

Thanks to this theorem the following formulas
$f(x)\leq g(y)$ and $f(x)=g(y)$ are satisfied if and only if the pair $(x,y)$ belongs respectively to $\{(x,y) \in \mathbb{R}^{n}\times\mathbb{R}^{n}\; |\; f(x) \leq g(y)\}$ and $\{(x,y) \in \mathbb{R}^{n}\times\mathbb{R}^{n}\; |\; f(x) =g(y)\}$. Therefore $f(x)=g(y)$ and $f(x)=g(y)$ can be treated as first-order formulas.

Moreover, notice that a structure on $(\mathbb{R},+,0)$ is compatible with addition and multiplication. This means that the sets $\{(x,y,z)\in \mathbb{R}^{3}\;|\; x+y=z\}$ and $\{(x,y,z)\in \mathbb{R}^{3}\;|\; x\cdot y=z\}$ are definable. Therefore the sum, the difference, the product and the quotient of definable functions is a definable function.

\begin{theorem}[\cite{long} Theorem 1.11]
\label{theo1}
The Euclidean distance from a semialgebraic set is semialgebraic. 
\end{theorem}

Notice that the $l_{1}$ norm and the supremum norm are both semialgebraic so the theorem above holds with the coresponding two metrics too.

\begin{definition}[\cite{multi}]
Let $A \subset \mathbb{R}^{n}$ and $f:A \longrightarrow \mathbb{R}$ be a function. For $f$ we define \textit{the generalized set of zeroes} as $f^{-1}(0)^{\gamma}:=\pi( \;\overline{\Gamma_{f}} \cap \mathbb{R}^{n}\times\{0\})$, where $\pi:\mathbb{R}^{n}\times \mathbb{R} \longrightarrow \mathbb{R}^{n}$ is the natural projection.
\end{definition}

If $f$ is continuous, the set of zeroes and the generalized set of zeroes are the same. If $f$ is a definable function then the generalized set of zeroes of $f$ is a definable set.
Moreover, we can also describe this set equivalently by:
\[f^{-1}(0)^{\gamma} \ni x_{0} \Longleftrightarrow \; \exists \; A \ni x_{\nu} \longrightarrow x_{0}: f(x_{\nu}) \longrightarrow 0 .  \]

\begin{remark}
\label{remark}
Due to the Curve Selection Lemma, in the definable case \\
$x_{0} \in f^{-1}(0)^{\gamma} \Longleftrightarrow \exists \; \gamma:(0,\epsilon) \longrightarrow A$ a definable curve such that $\lim_{t \to 0} \gamma(t)=x_{0}$ and $\lim_{t \to 0}f(\gamma(t))=0$.\\
Indeed, for $(x_{0},0) \in \overline{\Gamma_{f}}$ we either have $x_{0} \in A$ and $f(x_{0})=0$ in which case $\gamma(t) \equiv 0$ or one of two remaining possibilities: $x_{0} \in A$ with $f(x_{0}) \neq 0$ or $x_{0} \notin A$. In both these cases $(x_{0},0) \in \overline{\Gamma_{f}} \backslash \Gamma_{f}$ where by the Curve Selection Lemma there is a definable curve $\tilde{\gamma}:(0,\epsilon) \longrightarrow \Gamma_{f}$ with 
$ \tilde{\gamma}(t) \longrightarrow (x_{0},0)$ as $t \longrightarrow 0$. Clearly, for $\pi(x,y)=x$ we obtain a curve $\gamma(t)=\pi((\tilde{\gamma}(t))$ (since $\tilde{\gamma}$ cannot be projected into a singleton) satisfying the requirements.
\end{remark}

\begin{theorem}[\cite{long} Theorem 2.14](Łojasiewicz inequality)
\label{thm11}
Let $f,g:K \longrightarrow \mathbb{R}$ be two continuous semialgebraic functions on a compact set $K$ such that $f^{-1}(0) \subset g^{-1}(0)$.\\ Then there exist $C,N>0$ such that:
\[|f(x)|\geq C|g(x)|^{N}, \quad x \in K.\]
\end{theorem}

\begin{proposition}[\cite{long} Corollary 2.16]
\label{prop6}
Let $A,B$ be two non-empty compact semialgebraic sets. Then for some $N,C>0$ there holds:
\[dist(x,A) \geq C dist(x,A \cap B)^{N}, \quad x \in B.\]
\end{proposition}

\section{Generalized Łojasiewicz's inequality}
In \cite{categ} we can find version of Theorem \ref{thm11} for any definable structure. The aim of this chapter will be to prove that under some conditions we can omit the continuity assumption for the function $f$ and require from the function $g$ only to be bounded. Moreover, we will show that adding one condition we can formulate the theorem for any definable set not necessarily compact or even bounded.\bigskip \\
Let $\Phi_{p}$ be the set of all odd, increasing bijections $\phi:\mathbb{R} \longrightarrow \mathbb{R}$ that are of class $\mathcal{C}^{p}$ and $p$-flat at 0, meaning that $ \phi^{(k)}(0)=0$ for $k \in \{1,...,p\}$. According to the ideas from [5], in what follows $p\geq 1$ may be any fixed a priori integer number.\bigskip\\
To prove the announced results we will need the following theorems and propositions.

\begin{theorem}[\cite{multi} Theorem 2.3]
\label{zer}
If $f: \Omega \longrightarrow \mathbb{R}^{m}$, $\Omega \subset \mathbb{R}^{n}$ is a definable function, then for any point in the general set of zeroes of $f$ there is a neighbourhood $U$ and definable $\phi \in \Phi_{p}$ such that
\[||f(x)||\geq\phi(d(x,f^{-1}(0)^{\gamma})),\;\; x \in U \cap \Omega.\]
\end{theorem}

\begin{proposition}[\cite{categ} C.7 Lemma] \label{propw}
Let $f:A \times (0,\infty) \longrightarrow \mathbb{R}$, $A \subset \mathbb{R}^{n}$ be a definable function. Then there exists a definable function $\phi \in \Phi_{p}$ such that $\lim_{t \to 0^{+}}\phi(t)f(x,t)=0$ for each $x \in A$.
\end{proposition}

\begin{theorem}[\cite{categ} 4.20]
\label{nier} Let $f,g:A \longrightarrow \mathbb{R}$ be  continuous, definable functions on a compact set $A$ with $f^{-1}(0)\subset g^{-1}(0)$. Then there exists a definable $\phi \in \Phi_{p}$ such that $|\phi(g(x))|\leq |f(x)|$ for all $x \in A$.
\end{theorem}

\begin{proposition}
\label{minimum}
Let $\varphi_{1},...,\varphi_{n}: \mathbb{R} \longrightarrow \mathbb{R}$ be continuous, definable, odd functions. Define $\varphi(t):= \text{sgn}(t)\min{\{|\varphi_{1}(t)|,...,|\varphi_{n}(t)|\}}$, $t \in \mathbb{R}$. There exist $\epsilon >0$ and $i_{0} \in \{1,...,n\}$ such that $\varphi_{i_{0}}\arrowvert_{[- \epsilon,\epsilon]}=\varphi \arrowvert_{[- \epsilon,\epsilon]}$.
\begin{proof}
First note that the function $\varphi$ is continuous, odd and definable.
The proof will be by induction on $n$. For $n=1$ it is obvious. Let us consider the case $n=2$.  The set $\pi_{1}(\Gamma_{\varphi_{1}}\cap\Gamma_{\varphi_{2}})$ is not empty because $\varphi_{1}(0)=0=\varphi_{2}(0)$. Since the functions are definable, $\pi_{1}(\Gamma_{\varphi_{1}}\cap\Gamma_{\varphi_{2}})$ is a finite union of points and intervals. If there exists $c>0$ such that $[0,c] \subset \pi_{1}(\Gamma_{\varphi_{1}}\cap\Gamma_{\varphi_{2}})$, we have $\varphi\arrowvert_{[-c,c]}=\varphi_{1}\arrowvert_{[-c,c]}=\varphi_{2}\arrowvert_{[-c,c]}$ because the functions are odd. If not, denote $\epsilon:=\text{min}\{\left(\pi_{1}(\Gamma_{\varphi_{1}}\cap\Gamma_{\varphi_{2}})\backslash\{0\}\right) \cap \mathbb{R}_{+}\}$. Then one of the following cases holds: $\varphi_{1}>\varphi_{2}$ or $\varphi_{2}<\varphi_{1}$ on $(0,c]$. This means that there exists $i \in \{1,2\}$ such that $\varphi\arrowvert_{[0,c]}=\varphi_{i}\arrowvert_{[0,c]}$ and then $\varphi\arrowvert_{[-c,c]}=\varphi_{i}\arrowvert_{[-c,c]}$. Now take $n>0$ and suppose the theorem is true for all $m$ smaller than $n$. Define $\phi(t):= \text{sgn}(t)\min{\{|\varphi_{1}(t)|,...,|\varphi_{n-1}(t)|\}}$ for $t \in \mathbb{R}$. We can write $\varphi$ in the form $\varphi(t)= \text{sgn}(t)\min\{|\phi(t)|,|\varphi_{n}(t)|\}$. From the case $n=2$ we have two possibilities:\\
1) there exists $\epsilon>0$ such that $\varphi\arrowvert_{[-\epsilon,\epsilon]}=\varphi_{n}\arrowvert_{[-\epsilon,\epsilon]}$ and the theorem holds,\\
2) there exists $\epsilon_{1}>0$ such that $\varphi\arrowvert_{[-\epsilon_{1},\epsilon_{1}]}=\phi\arrowvert_{[-\epsilon_{1},\epsilon_{1}]}$.\\
Since $\phi(t)= \text{sgn}(t)\min{\{|\varphi_{1}(t)|,...,|\varphi_{n-1}(t)|\}}$, from the inductive hypothesis we know that there exist $\epsilon_{2}>0$ and $i_{0} \in \{1,...,n-1\}$ such that $\phi\arrowvert_{[-\epsilon_{2},\epsilon_{2}]}=\varphi_{i_{0}}\arrowvert_{[-\epsilon_{2},\epsilon_{2}]}$. Now taking $\epsilon = \min\{\epsilon_{1},\epsilon_{2}\}$ we get that $\varphi\arrowvert_{[-\epsilon,\epsilon]}=\varphi_{i_{0}}\arrowvert_{[-\epsilon,\epsilon]}$, which ends the proof of the theorem.
\end{proof}
\end{proposition}

Our following theorem is a general counterpart of \cite{multi} Theorem 2.3 and a definable counterpart of the semi-algebraic Lemma 2.3 from \cite{pre}.

\begin{theorem}
\label{uog}
Let $f,g:K \longrightarrow \mathbb{R}$ be two definable functions on a compact set $\emptyset \neq K \subset \mathbb{R}^{m}$. Moreover assume that $g$ is continuous and $f^{-1}(0)^{\gamma} \subset g^{-1}(0)$. There exists definable $\phi \in \Phi_{p}$ such that $|\phi(g(x))| \leq |f(x)|$ for all $x \in K$.
\begin{proof}
For any $\phi\in \Phi_p$ the function $h:=\phi(d(x,f^{-1}(0)^{\gamma}))$ is continuous and definable. Moreover $\phi(d(x,f^{-1}(0)^{\gamma}))=0 \Longleftrightarrow x \in f^{-1}(0)^{\gamma}$. For such a point $x$, from the assumptions we get that $g(x)=0$, so $h^{-1}(0) \subset g^{-1}(0)$. Fix now a point $a \in f^{-1}(0)^{\gamma}$ and take $U_{a}, \phi_{a}$ as in Theorem \ref{zer}. Notice that the function $d(\cdot,f^{-1}(0)^{\gamma})$ extends continuously to $\overline{U_{a}\cap K}$ and $h^{-1}(0) \subset g^{-1}(0)$ on the extension. The set $\overline{U_{a}\cap K}$ is compact, so from Theorem \ref{nier} we will find $\psi_{a} \in \Phi_{p}$ such that
\[d(x,f^{-1}(0)^{\gamma}) \geq |\psi_{a}(g(x))| \;\text{for}\; x \in \overline{U_{a}\cap K}.\]
We get that
\[|f(x)|\geq \phi_{a}(d(x,f^{-1}(0)^{\gamma})) \geq  \phi_{a}(|\psi_{a}(g(x))|) \;\text{for}\; x\in U_{a}\cap K,\; \phi_{a} \; \text{being increasing}.\]
Because $\phi_{a}$ is an odd function, we have \[\phi_{a}(|\psi_{a}(g(x))|)=|\phi_{a}(\psi_{a}(g(x)))|.\]
Taking $\varphi_{a}= \phi_{a} \circ \psi_{a}$ we obtain a function from $\Phi_{p}$ such that
\[|f(x)|\geq   |\varphi_{a}(g(x))| \;\text{for}\; x\in U_{a}\cap K.\]
Now notice that $f^{-1}(0)^{\gamma}=\overline{\Gamma_{f}}\cap K \times \{0\}$. The set $f^{-1}(0)^{\gamma}$ is a closed subset of $K$ thus compact. Let $\{U_{a}\}_{a \in f^{-1}(0)^{\gamma}}$ be an open cover of $f^{-1}(0)^{\gamma}$. Because $f^{-1}(0)^{\gamma}$ is compact we can find a finite subcover $\{U_{a_{1}},...,U_{a_{n}}\}$. From proposition \ref{minimum} we know that there exist $1>\epsilon >0$ and $i_{0} \in \{1,...,n\}$ such that $|\varphi_{a_{i_{0}}}| \leq |\varphi_{a_{i}}|$ on $[-\epsilon,\epsilon]$ for all $i \in \{1,...,n\}$. Since $g$ is continuous function on compact set there exists $t_{0}>0$ such that $|g(x)| \leq t_{0}$ on $K$. We may assume that $\epsilon<t_0$. Define the function 
\[\varphi(x) := \varphi_{a_{i_{0}}}(\frac{\epsilon}{t_{0}}x).\]
Take  $t \in \mathbb{R}$ such that $|t| \in [0,t_{0}]$. Notice that $|\frac{\epsilon}{t_{0}} t| \in [0,\epsilon]$ and $|\varphi_{a_{i}}(\frac{\epsilon}{t_{0}}t)| \leq |\varphi_{a_{i}}(t)|$ for all $\varphi_{a_{i}}$ because the functions are increasing and odd. We get that
\[|\varphi(t)| = |\varphi_{a_{i_{0}}}(\frac{\epsilon}{t_{0}}t)| \leq |\varphi_{a_{i}}(\frac{\epsilon}{t_{0}}t)| \leq |\varphi_{a_{i}}(t)|,\]
and then
\[|f(x)|\geq |\varphi(g(x))|,\; x \in \left(\bigcup_{i=1}^{n}U_{a_{i}}\right)\cap K. \]
Now consider the set $B:= K \backslash \bigcup_{i=1}^{n}U_{a_{i}}$. It is closed thus compact. We will prove that $|f|$ is bounded from below on $B$. Suppose to the contrary that there exists a sequence $\{x_{k}\}_{k=1}^{\infty} \subset B$ such that $|f(x_{k})| \longrightarrow 0$. Since $B$ is compact there exists a subsequence of $\{x_{k}\}_{k=1}^{\infty}$ tending to some $x_{0}$. This means that $x_{0} \in f^{-1}(0)^{\gamma} \subset g^{-1}(0)$. On the one hand $x_{0} \in B$ from the closedness of $B$, while on the other hand $x_{0} \in\bigcup_{i=1}^{n}U_{a_{i}}$ because $\bigcup_{i=1}^{n}U_{a_{i}}$ is open cover of $f^{-1}(0)^{\gamma}$. We get that 
\[x_{0} \in \left(\bigcup_{i=1}^{n}U_{a_{i}}\right)\cap B =\emptyset\]
which is a contradiction. Take $M>0$ such that $|f|>M$ on $B$. Define $L:=|\varphi(t_{0})|$. Take $C:=\min\{\frac{M}{L},1\}$ and define $\phi:=C\varphi$. We get the estimation
\[|f(x)|\geq |\phi(g(x))|\; \text{for}\; x\in K.\]
\end{proof}
\end{theorem}

Now we will present the most generalized version of Łojasiewicz's inequality.

\begin{theorem}
\label{glowne}
Let $f,g:A \longrightarrow \mathbb{R}$, $A \subset \mathbb{R}^{m}$ be two definable functions such that $|g|$ is bounded. Moreover assume the following condition:\\
\[ \forall \{x_{n}\}_{n=0}^{\infty} \subset A: f(x_{n}) \longrightarrow 0 \Longrightarrow g(x_{n}) \longrightarrow 0.\]
There exists definable $\varphi \in \Phi_{p}$ such that $|f(x)| \geq |\varphi(g(x))|$ for $x \in A$.
\begin{proof}
Define the following function:
\[\psi: (0, \infty) \ni t \longrightarrow sup\{ \frac{1}{|f(x)|}: x \in A \cap g^{-1}(t) \}.\]
Observe that for $t \in (0,\infty)$, $\frac{1}{|f(x)|}$ is well defined on $A \cap g^{-1}(t)$ because of $f^{-1}(0) \subset g^{-1}(0)$. The function $\psi$ is definable, it follows directly from the description of the graph, as the supremum in the definition is always finite. To prove this suppose there is a sequence $\{x_{n}\}_{n=1}^{\infty} \subset A \cap g^{-1}(t)$ such that $\frac{1}{|f(x_{n})|} \longrightarrow \infty$. We get that $f(x_{n}) \longrightarrow 0$. From the latter assumption we know that $g(x_{n}) \longrightarrow 0$ but $g(x_{n})=t$ for all $n \in \mathbb{N}$ so we get a contradiction.
Identifying $\psi$ with $\tilde{\psi}: \{0\} \times (0,\infty) \ni (0,t) \longrightarrow \psi(t) \in \mathbb{R}$, by proposition \ref{propw}, we can find a function $\phi \in \Phi_{p}$ such that $\lim_{t \to 0^{+}}\phi(t)\psi(t)=0$.
To prove the theorem it is sufficient to show that the function $\frac{|\phi \circ g|}{|f|}$ is bounded on $A \backslash \{g^{-1}(0)\}$. Suppose to the contrary that there is a sequence $\{x_{n}\}_{n=1}^{\infty} \subset A$ such that $\frac{|\phi(g(x_{n}))|}{|f(x_{n})|} \longrightarrow \infty$. The numerator is bounded because $g$ is bounded and $\phi$ is increasing. We get that $|f(x_{n})| \longrightarrow 0$. From the assumption we have that $g(x_{n}) \longrightarrow 0$. We get the following inequalities:
\[0 \leq \frac{|\phi(g(x_{n}))|}{|f(x_{n})|} \leq |\phi(g(x_{n}))\psi(g(x_{n}))| \longrightarrow 0.\]
It gives us a contradiction with the fact that $\frac{|\phi(g(x_{n}))|}{|f(x_{n})|} \longrightarrow \infty$.
This means that there exists $C>0$ such that
\[\frac{|\phi(g(x_{n}))|}{|f(x_{n})|} \leq C \]
and then 
\[\frac{1}{C} |\phi(g(x))| \leq |f(x)| \; \text{on} \; A \backslash g^{-1}(0).\]
Notice that this inequality is obviously satisfied for $x \in g^{-1}(0)$ so taking $\varphi:= \frac{1}{C}\phi$, we get that 
\[|\varphi(g(x))| \leq |f(x)| \; \text{for} \; x \in A.\]
\end{proof}
\end{theorem}

Observe that in the theorem above the condition $f^{-1}(0) \subset g^{-1}(0)$ is realised by constant sequences.

\begin{example}
The assumption of boundedness of the function $g$ cannot be omitted. Consider the following semialgebraic functions: $f: [0,1] \ni x \longrightarrow x \in \mathbb{R}$ and 
\[g(x) = \left\{ \begin{array}{ll}
\frac{1}{1-x}, & \textrm{when $x \in [0,1),$}\\
1, & \textrm{when $x=0$}.\\
\end{array} \right.\]
Suppose there exists $\phi \in \Phi_{p}$ such that $\phi(g(x)) \leq f(x),$ for $x \in [0,1]$. Take a sequence $\{x_{n}\}_{n=1}^{\infty} \subset [0,1)$ such that $x_{n} \longrightarrow 1$.
Notice that $f(x_{n})<1$ for all $n$. On the other hand $g(x_{n}) \longrightarrow \infty$ so $\phi(g(x_{n}))$ because $\phi$ is a continuous increasing bijection.
\end{example}

\begin{example}
Notice that we cannot ommit the assumption that $f^{-1}(0) \subset g^{-1}(0)$. In particular we cannot replace it just by the condition $f^{-1}(0)^{\gamma} \subset g^{-1}(0)^{\gamma}$. Let $f:[0,1] \longrightarrow 0 \in \mathbb{R}$ and
\[g(x) = \left\{ \begin{array}{ll}
0, & \textrm{when $x \in (0,1],$}\\
1, & \textrm{when $x=0$}.\\
\end{array} \right.\]
Observe that we have $g^{-1}(0)^{\gamma}=f^{-1}(0)^{\gamma}=[0,1]$. Suppose there exists $\phi \in \Phi_{p}$ such that $\phi(g(x)) \leq f(x),$ for $x \in [0,1]$. For $x=0$ we have that $g(0)=1$ and $\phi(g(0))=\phi(1) \leq f(0)=0$ which gives us a contradiction with the fact that $\phi$ is increasing odd bijection.
\end{example}

\section{Generalized Łojasiewicz's inequality for subanalytic functions}

The last thing of this part will be to formulate the generalized Łojasiewicz's inequality in the case of subanalytic sets. We have to consider this case separately because subanalytic sets do not from an o-minimal structure. For more information we refer the reader to [3].

\begin{definition}
A set $A\subset \mathbb{R}^{n}$ is called \textit{semi-analytic}, if for any $ x \in \mathbb{R}^{n}$, there are a neighbourhood $U \ni x$ and analytic functions $f_{i}, g_{ij}$ in $U$ such that
\[A \cap U = \bigcup_{i=1}^{p}\bigcap_{j=1}^{q}\{x \in U | f_{i}(x)=0, g_{ij}(x) >0\}.\]
We can define this in a more general setting taking an analytic variety $M$ instead of $\mathbb{R}^{n}$.

\end{definition}

\begin{definition}
A set $E \subset M$, $M$ real analytic variety, is called \textit{subanalytic} if for any $x \in M$ there is a neighbourhood $U \ni x$ such that $E \cap U = \pi(A)$, where $\pi:M \times N \longrightarrow M$ is the natural projection, $N$ is a real variety and $A$ is semi-analytic and relatively compact in $M \times N$.
\end{definition}

\begin{definition}
A map $f:E \longrightarrow N$, where $E\subset M$ is a nonempty subanalytic set, is subanalytic if its graph is subanalytic in $M \times N$
\end{definition}

Basic operations on subanalytic sets like taking union, intersection, complement and closure preserves subanalyticity \cite{long}.

\begin{theorem}[\cite{loj} Propoition 3.17]
\label{thm13}
Let $f,g:K \longrightarrow \mathbb{R}$ be two continuous subanalytic functions on a compact set $K$ such that $f^{-1}(0) \subset g^{-1}(0)$.\\ Then there exist $C,N>0$ such that:
\[|f(x)|\geq C|g(x)|^{N}, \quad x \in K.\]
\end{theorem}

\begin{theorem}[\cite{multi} Theorem 2.1]
\label{thm14}
If $f: \Omega \longrightarrow \mathbb{R}$ is a subanalytic function, $\Omega \subset \mathbb{R}^{n}$ then for any point a in the general set of zeroes of $f$ there is a neighbourhood $U$ and positive constants $C \; \text{and} \; \alpha$,  such that
\[|f(x)|\geq Cd(x,f^{-1}(0)^{\gamma})^{\alpha},\;\; x \in U \cap \Omega.\]
\end{theorem}

\begin{theorem}
\label{subuog}
Let $f,g:K \longrightarrow \mathbb{R}$ be two subanalytic functions on a compact set $\emptyset \neq K \subset \mathbb{R}^{n}$. Moreover assume that $g$ is continuous and $f^{-1}(0)^{\gamma} \subset g^{-1}(0)$. There exist $C,\alpha>0$ such that $C|g(x)|^{\alpha} \leq |f(x)|$ for all $x \in K$.
\begin{proof}
At the beginning, using theorems \ref{thm13} and \ref{thm14} we argue in the same way as in the proof of theorem \ref{uog}. We get a finite cover $\{U_{1},...,U_{n}\}$ of $f^{-1}(0)^{\gamma}$ and positive constants $C_{1},...,C_{n},\alpha_{1},...,\alpha_{n}$ such that
\[C_{i}|g(x)|^{\alpha_{i}} \leq f(x),\; x \in U_{i}.\]
Moreover we know that $f$ is bounded from below on $K \backslash \left(\bigcup_{i=1}^{n} U_{i} \right)$ by some constant $D$. Define $\alpha = max \{\alpha_{1},...,\alpha_{n},1\}$. Notice now that if $|g(x)| \leq 1$ for some $x \in K$ then $|g(x)|^{\alpha} \leq |g(x)|^{\alpha_{i}}$ for $i=1,...,n$. Since $|g|^{\alpha}$ is continuous on compact set there exists $C_{0}>0$ such that $|g(x)|^{\alpha} \leq C_{0}$ for $x \in K$. Defining $C= min\{\frac{D}{C_{0}},C',\frac{C'}{C_{0}}\}$, where $C'= min\{C_{1},...,C_{n}\}$. We get
\[C|g(x)|^{\alpha} \leq |f(x)|, \; x \in K.\]
\end{proof}
\end{theorem}

\begin{theorem}
\label{subglowne}
Let $f,g:A \longrightarrow \mathbb{R}$ be two subanalytic functions on compact subanalytic subset $A \subset \mathbb{R}^{m}$ such that $|g|$ is bounded. Moreover assume the following condition:\\
\[ \forall \{x_{n}\}_{n=0}^{\infty} \subset A: f(x_{n}) \longrightarrow 0 \Longrightarrow g(x_{n}) \longrightarrow 0. \;\;\; (*)\]
Then there exist constants $C, \alpha >0$ such that $|f(x)| \geq C|g(x)|^{\alpha}$ for $x \in A$.
\begin{proof}
In this proof we will consider the norm in $\mathbb{R}^{2}$ given by $||(x,y)||=|x|+|y|$. Following \cite{loj} (Prop 3.17), define the following sets
\[B = \{(s,t) \in  \mathbb{R} \; | \; \exists x\in A:\; f(x)=s, \; g(x)=t\},\]
\[B' = \{(x,s,t) \in A \times \mathbb{R}^{2} \; |\; f(x)=s,\; g(x)=t\}. \]
Let $\pi:\mathbb{R}^{3} \longrightarrow \mathbb{R}^{2}$ be the projection onto the last two coordinates. We have $\pi(B')=B$ and $\pi$ is bounded in the direction of the projection, so $B$ is subanalytic \cite{loj} (Prop 2.9).\\
Suppose now that $f^{-1}(0)^{\gamma} \neq \emptyset$.
Notice that $\overline{B} \cap (\{0\}\times \mathbb{R}) = \{ (0,0)\}$. Indeed, take a sequence $\{(x_{n},y_{y})\}_{n=1}^{\infty}\subset B$ and suppose that $ x_{n} \longrightarrow 0, \; y_{n} \longrightarrow y_{0}$ with $y_{0} \neq 0$. We get a contradiction with the condition ($*$).
From the regular separation theorem \cite{loj} (Thm 3.5) we get that there is a neighbourhood $U\subset \mathbb{R}^{2}$ of $(0,0)$ and constants $C_{1},\alpha>0$ such that
\[d((s,t),\{0\}\times\mathbb{R}) \geq C_{1}d((s,t),\overline{B}\cap(\{0\}\times \mathbb{R}))^{\alpha}=C_{1}d((s,t),\{(0,0)\})^{\alpha},\; \text{for}\; (s,t) \in U.\]
Define now the sets $Z_{1}=\pi(\Gamma_{(f,g)})\cap U$ and $Z_{2}=\pi(\Gamma_{(f,g)})\backslash U$. Obviously we have $\pi(\Gamma_{(f,g)}) = Z_{1} \cup Z_{2}$. Notice that the first coordinate in the set $Z_{2}$ is separated from $0$. To see this suppose it is not the case and take a sequence $\{(x_{m},y_{m})\}_{m=1}^{\infty} \subset Z_{2}$ such that $x_{m} \longrightarrow 0$. From the condition ($*$) we get that $y_{n} \longrightarrow 0$. We have that $\{(x_{m},y_{m})\}_{m=1}^{\infty} \subset \mathbb{R}^{2} \backslash U$, so we get a contradiction because $(0,0) \in \text{int}\; U$.
Take $K>0$ such that for every $(s,t) \in Z_{2}$ we have $s>K$. We know that $|g|$ is bounded on $A$ whence so is $|g|^{\alpha}$. Take $L>0$ such that $L>|g(x)|^{\alpha}$, for every $x \in A$. We can find $C_{2}$ such that $K>C_{2}L$. Take $C=\min\{C_{1},C_{2}\}$.\\
Notice that for $ x \in A$ such that $(f(x),g(x))=(s,t) \in Z_{1}\subset U$, we have
\[|f(x)|=d((s,t),\{0\}\times\mathbb{R}) \geq C_{1}d((s,t),\{(0,0)\})^{\alpha}=C_{1}(|s|+|t|)^{\alpha} \geq C_{1}|t|^{\alpha} \geq C|g(x)|^{\alpha}.\]
On the other hand, for $ x \in A$ such that $(f(x),g(x))=(s,t) \in Z_{2}\subset U$ we have
\[|f(x)| \geq K \geq C_{2}L \geq CL \geq C|g(x)|^{\alpha}.\]
Since for $x \in A$, there holds $(f(x),g(x)) \in Z_{1}$ or $(f(x),g(x)) \in Z_{2}$, this part of the proof is done.
If $f^{-1}(0)^{\gamma}=\emptyset$ then $|f|$ is bounded from below and we are reasoning as before. It is bounded from below because otherwise we would have a sequence of arguments on which the function $f$ would tend to zero, and since $A$ is compact, the generalized set of zeroes could not be empty.  
\end{proof}
\end{theorem}

\begin{theorem}
\label{subuog2}
Let $f,g:A \longrightarrow \mathbb{R}$ be two subanalytic functions on bounded subanalytic subset $A \subset \mathbb{R}^{m}$ such that $|g|$ is bounded. Moreover, assume the following condition:\\
\[ \forall \{x_{n}\}_{n=0}^{\infty} \subset A: f(x_{n}) \longrightarrow 0 \Longrightarrow g(x_{n}) \longrightarrow 0. \;\;\; (*)\]
Then there exist constants $C, \alpha >0$ such that $|f(x)| \geq C|g(x)|^{\alpha}$, for $x \in A$.
\begin{proof}
First notice that $\overline{A}$ is compact and subanalytic.
We define the following subanalytic functions:
\[\widetilde{g}(x) = \left\{ \begin{array}{ll}
0, & \textrm{when $x \in \overline{A} \backslash A,$}\\
g(x), & \textrm{when $x \in A$}.\\
\end{array} \right.\]
\[\widetilde{f}(x) = \left\{ \begin{array}{ll}
0, & \textrm{when $x \in \overline{A} \backslash A,$}\\
f(x), & \textrm{when $x \in A$}.\\
\end{array} \right.\]
The functions $\widetilde{g},\widetilde{f}$ satisfy the assumptions of the previous theorem so we get the statement.
\end{proof}
\end{theorem}

There is a natural question whether we can have a counterpart of the theorem \ref{glowne} in the subanalytic case. The answer turns out to be negative, as shown in the example which has been developed in collaboration with Maciej Denkowski.
\begin{example}
\label{ex}
Define the following subanalytic functions:
\[g(x)=x-\floor*{x}, \; \text{for} \; x \in [0,\infty)\]
and
\[f(x)=(x-\floor*{x})^{\floor*{x}}, \; \text{for} \; x \in [0,\infty). \]
Notice that the graph of the function $g$, restricted to the interval $[n,n+1)$ for $n \in \mathbb{N}$, is just translation of the graph of the function
\[\tilde{g}(x)=x, \; \text{for} \; x \in [0,1).\]
Similarly, the graph of the function $f$, restricted to the interval $[n,n+1)$ for $n \in \mathbb{N}$, is just translation of the graph of the function
\[\tilde{f}_{n}(x)=x^{n}, \; \text{for} \; x\in[0,1).\]
Supposing that there is a counterpart of theorem \ref{glowne} in the subanalytic case, we will have $C, \alpha>0$ such that
\[f(x) \geq Cg(x)^{\alpha}, \; \text{for} \; x \in [0,\infty).\] This means that the inequality
\[ x^{n} \geq Cx^{\alpha}, \; \text{for} \; x \in [0,1)\] would be satisfied for all $n \in \mathbb{N}$, what yields a contradiction.
\end{example}

A result similar to Theorem \ref{subuog2} is to be found in \cite{valette}:

\begin{proposition}[\cite{valette} Prop 1.1]
\label{val}
Let $f$ and $g$ be two globally subanalytic functions on a globally subanalytic set $A$ with $\sup_{x \in a} |f(x)| < \infty$. Assume that $\lim_{t \to 0} f(\gamma(t))=0$ fro every globally subanalytic arc $\gamma:(0,\epsilon) \longrightarrow A$ satisfying $\lim_{t \to 0} g(\gamma(t))=0$.
Then there exist $N \in \mathbb{N}$ and $C \in \mathbb{R}$ such that for any $x \in A$:
\[ |f(x)|^{N} \leq C|g(x)|.\]
\end{proposition}

The proof extends easily to polynomially bounded o-minimal structures. There are two major differences with Theorem \ref{subuog2}. In Proposition \ref{val} we do not require $A$ to be bounded but instead both the functions have to be globally subanalytic. In Theorem \ref{subuog2} $g$ is in fact globally subanalytic, while $f$ need not be such. Example \ref{ex} illustrates this difference between the two results. The second difference is only in notation as is seen from the lemma below.

\begin{lemma}
\label{v1}
Let $f,g: A \longrightarrow \mathbb{R}$ be two sybanalytic functions on a bounded subanalytic set $A \subset \mathbb{R}^{n}$.
Then the following assertions are equivalent:\\
i) for any analytic $\gamma:(0,\epsilon) \longrightarrow A$ with $\lim_{t \to 0}f(\gamma(t))$ there is $\lim_{t \to 0} g(\gamma(t))=0$,\\
ii) for any $(x_{\nu})\subset A$ such that $f(x_{\nu}) \longrightarrow 0 $ there is $ g(x_{\nu}) \longrightarrow 0$.
\begin{proof}
We only need to prove $(i) \Longrightarrow (ii)$. Suppose on the contrary that there exists a sequence $(x_{\nu}) \subset A$ such that $ |g(x_{\nu})|>c$ while $f(x_{\nu}) \longrightarrow 0$ where $c>0$ is a constant. Observe that $|g|$ is subanalytic (since $\Gamma_{|g|}=\left(\Gamma_{g} \cap A \times [0,+\infty)\right) \cup  \varphi \left(\Gamma_{g} \cap A \times [-\infty,0) \right)$ where $\varphi(x,y)=(x,-y)$ is bi-analytic) whence the set $A_{c}=\{x \in A \; | \; |g(x)|>c\}$ is subanalytic too, being the complement of $\{x \in A \; | \; |g(x)| <x\}=\pi(\Gamma_{|g|} \cap A \times[0,c))$ where $\pi(x,y)=x$  on the projected subanalytic set has relativitely compact preimages in the $y$ direction. We do not need $|g|$ to be bounded.
The set $A$ being bounded, we may assume $ x_{\nu} \longrightarrow x_{0}$. Then $(x_{0},0) \in \overline{(A_{c}\times \{0\}) \cap \Gamma_{f}}$ and we may proceed exactly as in the Remark \ref{remark} keeping in mind that in the subanalytic case the Curve Selection Lemma gives an analytic curve.
\end{proof}
\end{lemma}

The lemma above can be carried over to the general definable case.

\begin{lemma}
\label{v2}
Let $f,g: A \longrightarrow \mathbb{R}$ be two definable functions.
Then the following assertions are equivalent:\\
i) for any definable $\gamma:(0,\epsilon) \longrightarrow A$ with $\lim_{t \to 0}f(\gamma(t))$ there is $\lim_{t \to 0} g(\gamma(t))=0$,\\
ii) for any $(x_{\nu})\subset A$ such that $f(x_{\nu}) \longrightarrow 0 $ there is $ g(x_{\nu}) \longrightarrow 0$.
\begin{proof}
We prove $(ii) \Longrightarrow (i)$. Note that $A$ may be unbounded. Therefore either $(x_{\nu})$ has a limit in $\overline{A}$ (after passing to a subsequence and then we may repeat the preceding prove based on the Curve Selection Lemma, or $||x_{\nu}|| \longrightarrow +\infty$. Then we remark that $\varphi(x)= \frac{x}{||x||^{2}}$ is a semialgebraic diffeomorphism $\mathbb{R}^{n} \backslash \{0\} \longrightarrow \mathbb{R}^{n} \backslash \{0\}$. Thus $\tilde{A}=\varphi(A\backslash \{0\})$ is definable and $\varphi(x_{\nu}) \longrightarrow 0$. Taking now the definable curve $\tilde{\gamma}:(0,\epsilon) \longrightarrow \tilde{A}$ satisfying $\lim_{t \to 0} \tilde{\gamma}(t)=0$, we obtain a definable curve $\gamma=(\varphi^{-1}) \circ \tilde{\gamma}:(0,\epsilon) \longrightarrow A$, for which $\lim_{t \to 0} \gamma(t)=+\infty$. This allows us to apply the first part of the proof to $\tilde{A}$ and for $f \circ \varphi^{-1}$, $g \circ \varphi^{-1}$ which ends the argument
\end{proof}
\end{lemma}

\begin{remark}
In completion to A.Valette's result, we note that in Lemma \ref{v1} analytic curves are sufficient (but we keep the assumption that $A$ is bounded), while Lemma \ref{v2} is a direct complement to A.Vallete's definable version. 
\end{remark}

\section{Łojasiewicz's inequality for multifunctions}
Our next goal will be to use the previously proven theorem to formulate some counterparts of Łojasiewicz's inequality for definable multifunctions. The definitions presented  come from \cite{multi}.

\subsection{Definable multifunctions}

\begin{definition}
A \textit{multifunction} is a function $F: \mathbb{R}^{m} \longrightarrow \mathcal{P}(\mathbb{R}^{n})$. The elements of $F(a)$ are called the values of the multifunction $F$ at $x$. Moreover, we consider the domain of a multifunction as the set $dom \; F=\{x \in \mathbb{R}^{n} \; | \; F(x) \neq \emptyset\}.$\\
A multifunction whose graph $\Gamma_{F}=\{(x,y)\in \mathbb{R}^{m}\times \mathbb{R}^{n} \;| y \in F(x) \}$ is a definable set, is said to be a \textit{definable multifunction}.\\
The domain of a definable multifunction is a definable set.
\end{definition}

\begin{definition}
For a multifunction $F: \mathbb{R}^{m} \longrightarrow \mathcal{P}(\mathbb{R}^{n})$ and $a \in \overline{dom F\backslash \{a\}}$ we define the Kuratowski partial limits :
\begin{align*}
&\bullet \; y \in \liminf_{x \to a} F(x) \Leftrightarrow \forall \;dom\; F\backslash\{a\}\ni x_{\nu} \rightarrow a,\; \exists\; F(x_{\nu}) \ni y_{\nu} \rightarrow y \Longleftrightarrow \\
& \Leftrightarrow \ \forall U \ni y \ \text{neighborhood}, \exists V \ni a \ \text{neighborhood}, \forall x \in V \cap dom F\backslash \{a\}, U \cap  F(x) \neq \emptyset;\\
& \bullet \; y \in \limsup_{x \to a} F(x) \Leftrightarrow    \exists \;dom \;F\backslash\{a\}\ni x_{\nu} \rightarrow a,\; \exists\; F(x_{\nu}) \ni y_{\nu} \rightarrow y \Leftrightarrow \\
& \Leftrightarrow \;\forall \; U \ni y\; \text{neighborhood}, \forall \; V \ni a \; \text{neighborhood}, \exists \; x \in V\backslash \{a\}: \; U \cap F(x) \neq \emptyset.
\end{align*}
The inclusion $\liminf_{x \to a} F(x)\subset\limsup_{x \to a}F(x)$ always holds. If the reverse inclusion holds too, we have
$\limsup_{x \to a} F(x)=\liminf_{x \to a} F(x)$ \;and we talk about convergence in the Kuratowski sense.
\end{definition}

\begin{definition}
Let $a \in dom \;F$. With the previous notation we say that $F$ is:\\
i) outer semi-continuous in $a$, when $\limsup_{x \to a}F(x) \subset F(a)$;\\
ii) inner semi-continuous in $a$, when $\liminf_{x \to a}F(x) \supset F(a)$;\\
iii) upper semi-continuous in $a$, when $\limsup_{x \to a}F(x) = F(a)$;\\
iv) lower semi-continuous in $a$, when $\liminf_{x \to a}F(x) =F(a)$;\\
v) continuous in $a$, when $\liminf_{x \to a}F(x)=\limsup_{x \to a}F(x) =F(a)$.
\end{definition}

\begin{definition}[\cite{multi}]
For a multifunction $F:\mathbb{R}^{m} \longrightarrow P(\mathbb{R}^{n})$ and $a \in dom\;F$ we define:\\
- $F^{-1}(F(a)):=\{x \in \mathbb{R}^{m}\;|\;F(x)=F(a)\}$ -- the strong preimage,\\
- $F^{*}(F(a)):=\{x \in dom\;F\;|\;F(x)\subset F(a) \}$ -- the lower preimage,\\
- $F_{*}(F(a)):=\{x \in \mathbb{R}^{m}\;|\;F(a)\subset F(x)\}$ -- the upper preimage,\\
- $F^{\#}(F(a)):=\{x \in \mathbb{R}^{m} \;|\; F(x) \cap F(a) \neq \emptyset \}$ -- the weak preimage.
\end{definition}

\begin{proposition}[\cite{multi} Proposition 3.1]
For a definable multifunction $F: \mathbb{R}^{m} \longrightarrow P(\mathbb{R}^{n})$ and $a \in dom\;F$, all the types of defined above preimages are definable sets.
\end{proposition}

\subsection{Hausdorff distance}

\begin{definition}
Let $A,B$ be two non-empty, compact subsets of a metric space $X$. We define the Hausdorff metric as
\[dist_{H}(A,B):=\max\{\max_{x\in A}d(x,B),\max_{x\in B}d(x,A)\}.\]
\end{definition}

If the space $X$ is compact, then we can also consider the extended Hausdorff metric. We define then $dist_{H}(A,B)=diam\; X +1$, if $A$ or $B$ is the empty set but not both.

\begin{proposition}[\cite{multi} Lemma 5.1]
\label{prop4}
Let $f,g:A\longrightarrow \mathbb{R}$, $A \subset \mathbb{R}^{n}$ be two definable functions. Then the function $\phi(x)=max\{f(x),g(x)\}$ is definable.
\end{proposition}

\begin{proposition}[\cite{multi} Lemma 5.2]
\label{prop2}
Let $F: \mathbb{R}^{m} \longrightarrow \mathcal{P}(\mathbb{R}^{n})$ be a definable multifunction and $f:\mathbb{R}^{k}\times \mathbb{R}^{n} \longrightarrow \mathbb{R}$ a definable function. Then the function
\[\psi:\mathbb{R}^{k}\times dom\;F \ni(x,x') \longrightarrow \sup_{y \in F(x')} f(x,y) \in \mathbb{R} \cup \{\infty\}\]
is definable.
\end{proposition}

\begin{proposition}
\label{prop3}
Let $F: \mathbb{R}^{m} \longrightarrow \mathcal{P}(\mathbb{R}^{n})$ be a definable multifunction and $f:\mathbb{R}^{k}\times \mathbb{R}^{n} \longrightarrow \mathbb{R}$ a definable function. Then the function
\[\psi:\mathbb{R}^{k}\times domF \ni(x,x') \longrightarrow \inf_{y \in F(x')} f(x,y) \in \mathbb{R} \cup \{-\infty\}\]
is definable.
\begin{proof}
The infimum can be expressed with the supremum, so using \ref{prop2} we get the assertion.
\end{proof}
\end{proposition}

\begin{proposition}[\cite{multi} Proposition 5.3]
\label{prop7}
Let $F:\mathbb{R}^{m} \longrightarrow \mathcal{P}(\mathbb{R}^{n}),\;G:\mathbb{R}^{k} \longrightarrow \mathcal{P}(\mathbb{R}^{n})$ be two definable functions with compacts values. Then the function
\[\delta_{F,G}:dom\;F \times dom\;G \ni (x,x') \longrightarrow \max_{y\in G(x')}d(y,F(x)) \in \mathbb{R}\]
is definable.
\end{proposition}

\begin{proposition}[\cite{multi} Theorem 5.4]
\label{prop25}
Let $F:\mathbb{R}^{m} \longrightarrow \mathcal{P}(\mathbb{R}^{n}),G:\mathbb{R}^{k} \longrightarrow \mathcal{P}(\mathbb{R}^{n})$ be two definable multifuctions with compact values. Then the function:
\[d_{H}(F,G):dom\;F \times dom\;G \ni (x,x') \longrightarrow dist_{H}(F(x),G(x')) \in \mathbb{R}\]
is definable.
\end{proposition}

\begin{proposition}
\label{prop12}
Let $F:\mathbb{R}^{m} \longrightarrow \mathcal{P}(\mathbb{R}^{n})$ be a multifunction with compact values. Let $a\in dom\;F$ be a fixed point and $\{x_{\nu}\}_{\nu=1}^{\infty}\subset \mathbb{R}^{m}$ a sequence such that $x_{\nu} \longrightarrow x_{0}$ for some $x_{0} \in \mathbb{R}^{m}$. Then
\[dist_{H}(F(x_{\nu}),F(a)) \longrightarrow 0 \Longrightarrow \limsup_{x \to x_{0}}F(x) \supset F(a).\]
\begin{proof}[Proof]
Let $z \in F(a)$. Due to the fact that $dist_{H}(F(x_{\nu}),F(a))\longrightarrow 0$ holds, for every $\epsilon >0$ we have $dist_{H}(F(x_{k}),F(a))<\epsilon $ for sufficient large $k$. Then for fixed $k$ and $z \in F(a)$ we get that $d(z,F(x_{k})) \leq \max_{w \in F(a)}d(w,F(x_{k}))\leq \epsilon$ and the distance is realised for some $z_{k} \in F(x_{k})$ so we get a sequence $F(x_{k}) \ni z_{k} \longrightarrow z$. This means that $z \in \limsup_{x \to x_{0}}F(x)$.\end{proof}
\end{proposition}

\subsection{Łojasiewicz's inequality for multifunctions}

This subsection completes the results presented in \cite{multi}. In particular, it eliminates a flaw form

\begin{theorem}
\label{thm2}
Let $F:\mathbb{R}^{m} \longrightarrow \mathcal{P}(\mathbb{R}^{n})$ be an outer semi-continuous definable multifunction with compact values. For $a \in dom \;F$ there exist a definable $\varphi \in \Phi_{p}$ such that
\[dist_{H}(F(x),F(a)) \geq \varphi(d(x,F_{*}(F(a)))),\; x \in K\]
where $K$ is compact neighbourhood of $a$ such that $K$ in $dom \;F$.
\begin{proof}[Proof]
Let us define the functions $g(x):= d(x,F_{*}(F(a)))$, $f(x):=dist_{H}(F(x),F(a))$. Let $x_{\nu} \longrightarrow x_{0}$ be a sequence that $f(x_{\nu}) \longrightarrow 0$. We get that $\limsup_{x \to x_{0}}F(x) \supset F(a)$ by proposition \ref{prop12}. Due to the outer semi-continuity we have $\limsup_{x \to x_{0}}F(x) \subset F(x_{0})$. Thus there is an inclusion $F(a) \subset F(x_{0})$ and this means that $g(x_{0})=0$. Therefore, the functions $f,g$ satisfy the assumptions of the Theorem \ref{uog}, so we get the result sought after.
\end{proof}
\end{theorem}

\begin{example}
Our next goal will be to modify the assumptions of the last theorem in such a way that the thesis will hold for other types of preimages. Let us consider the following semialgebraic multifunction.\\

\[F(x) = \left\{ \begin{array}{ll}
\{\sqrt{x}\}, & \textrm{when $x \in (0,4)$}\\
\{2,1\}, & \textrm{when $x=4$}\\
\{-x+5\}, & \textrm{when $x \in (4,5).$}\\
\end{array} \right. \]
It is an upper semi-continuous semialgebraic multifunction, but for $a=1$ and the sequence $x_{n}=4+\frac{1}{n}$ we get $dist_{H}(F(x_{n}),F(a))\longrightarrow 0$, while $d(4,F^{-1}(F(1)))=d(4,\{1\})\neq 0$ and $d(4,F^{*}(F(1)))=d(4,\{1\})\neq 0$. This means that for the lower and strong preimages the assumption of upper semi-continuity is not enough. 
\end{example}

\begin{theorem}
Let $F:\mathbb{R}^{m} \longrightarrow \mathcal{P}(\mathbb{R}^{n})$ be an inner semi-continuous definable multifunction with compact values. For $a \in dom\; F$ there exist a definable $\varphi \in \Phi_{p}$ such that
\[dist_{H}(F(x),F(a)) \geq \varphi(d(x,F^{*}(F(a)))),\; x \in K\]
where $K$ is a compact neighbourhoood of $a$ such that $K$ in $dom \;F$.
\begin{proof}[Proof]
Define $g(x):= d(x,F^{*}(F(a)))$, $f(x):=dist_{H}(F(x),F(a))$. Let $x_{\nu} \longrightarrow x_{0}$ be a sequence such that $f(x_{\nu}) \longrightarrow 0$. We want to show that $F(x_{0}) \subset F(a)$. Fix $y \in F(x_{0})$. From the inner semi-continuity we get that for $x_{\nu} \longrightarrow x_{0}$ there exist a sequence $F(x_{\nu})\ni y_{\nu}$ such that $y_{\nu} \longrightarrow y$. From the fact that $dist_{H}(F(a),F(x_{\nu})) \longrightarrow~0$, and $F(a)$ is compact hence closed, we get that $y \in F(a)$. Therefore, the functions $f,g$ satisfy the assumptions of the Theorem \ref{uog} which give us the thesis.
\end{proof}
\end{theorem}

\begin{example}
Let us consider the semialgebraic multifunction, such that
\[G(x) = \left\{ \begin{array}{ll}
\{(x-2)^{2},1,2\}, & \textrm{when $x \in [-1,4)$}\\
\{1,2\}, & \textrm{when $x \in [4,6]$}.\\
\end{array} \right.\]
It is lower semi-continuous. Notice that for the sequence $x_{n}=4-\frac{1}{n}$ there is $dist_{H}(G(x_{n}),G(0)) \longrightarrow0$ but $d(4,G^{-1}(G(0))) \neq 0$. This means that the assumption of lower semi-continuity is not enough for the strong preimage. Previously we have showed that the upper-semicontinuity is not enough too. This leads to the following theorem. 
\end{example}

\begin{theorem}
Let $F:\mathbb{R}^{m} \longrightarrow \mathcal{P}(\mathbb{R}^{n})$ be a continuous definable multifunction with compact values. For $a \in dom\; F$ there exist a definable $\varphi \in \Phi_{p}$ such that
\[dist_{H}(F(x),F(a)) \geq \varphi(d(x,F^{-1}(F(a)))),\; x \in K\]
where $K$ is compact neighbourhood of $a$ such that $K \subset dom \;F$.
\begin{proof}[Proof]
Let us define the following two functions $g(x):= d(x,F^{-1}(F(a)))$ and $f(x):=dist_{H}(F(x),F(a))$. Let $x_{\nu} \longrightarrow x_{0}$ be a sequence such that $f(x_{\nu}) \longrightarrow 0$. We get that $\limsup_{x \to x_{0}}F(x) \supset F(a)$. Due to the continuity we have that $\limsup_{x \to x_{0}}F(x)=\liminf_{x \to x_{0}}F(x)= F(x_{0})$, so the inclusion $F(a) \subset F(x_{0})$ holds. A continuous multifunction is in particular inner semi-continuous, so from the proof of the previous theorem we know that the inclusion $F(x_{0})\subset F(a)$ holds too. We got the equality $F(x_{0})=F(a)$ which means that the functions $f$ and $g$ satisfy the assumptions of Theorem \ref{uog} what gives us the thesis.
\end{proof}
\end{theorem}

\begin{remark}
Notice that by using Theorem \ref{zer} and repeating the above proofs we can get versions of these three theorems for $x \in U \cap dom\;F$ where $U$ is some open neighbourhood of $a$. The difference is that $K$ can be any compact set that is contained in $dom\; F$ in particular not depending on $a$ while $U$ depends on $a$.
\end{remark}

There is a natural question whether we can use Theorem \ref{glowne} to get the previous theorems for any subset $K$, not necessarily compact. The next example provides the answer.

\begin{example}
Let us consider the following semialgebraic multifunction:
\[H(x) = \{x^{2}, x^{2}+1\}, \; \text{for} \; x \in (-1,1].\]
Taking $K$ as $(-1,1]$ and a sequence $x_{n}=-1+\frac{1}{n}$ we have $dist_{H}(H(x_{n}),H(1)) \longrightarrow 0$ but $H^{-1}(H(1))=H^{*}(H(1))=H_{*}(H(1))=\{1\}$ so distance from $x_{n}$ to these sets is not tending to $0$. The assumptions of Theorem \ref{glowne} are not satisfied.
\end{example}

\subsection{Multifunctions with closed values}

Our next goal is to generalize the presented facts to multifunctions with closed values. Let $\mathcal{F}_{n}$ be the family of all closed subsets of $\mathbb{R}^{n}$.

\begin{lemma}[\cite{multi} Lemma 6.1] 
Let $F:\mathbb{R}^{m} \longrightarrow \mathcal{P}(\mathbb{R}^{n})$ be a definable multifunction, and $h:\mathbb{R}^{n} \longrightarrow \mathbb{R}^{p}$ a definable function. Then the function $G: dom\;F\ni x \longrightarrow h(F(x)) \in \mathcal{P}(\mathbb{R}^{p})$ is definable.
\end{lemma}

\begin{lemma}[\cite{multi} Lemma 6.2]
Let $F,G: \mathbb{R}^{m} \longrightarrow \mathcal{P}(\mathbb{R}^{n})$ be two definable multifunctions. Then the functions $H_{1}=F(x)\cap G(x)$, $H_{2}=F(x)\cup G(x)$, $H_{3}=F(x)\backslash G(x)$ are definable.
\end{lemma}

\begin{theorem}[\cite{multi} Theorem 6.3]
There exist a metric $\text{dist}_{K}$ on the set $\mathcal{F}_{n}$ extending the Hausdorff metric such that for any two definable functions $F:\mathbb{R}^{m} \longrightarrow \mathcal{P}(\mathbb{R}^{n}), G:\mathbb{R}^{k}\longrightarrow \mathcal{P}(\mathbb{R}^{n})$ with compact values the induced function:
\[d_{K}(F,G):\mathbb{R}^{m}\times\mathbb{R}^{k} \ni (x,y) \longrightarrow dist_{K}(F(x),G(y)) \in \mathbb{R}\]
is definable.

\begin{proof}[Proof]
We identify $\mathbb{R}^{n}$ with $R:=\mathbb{R}^{n}\times\{-1\} \subset \mathbb{R}^{n+1}$ and consider the unit sphere $\mathbb{S}^{n}$ with the stereographic projection from the point on the pole $p=(0,...,0,1)$ to $R$,
\[s(x)=\left(\frac{-2 x_{1}}{x_{n+1}-1},...,\frac{-2 x_{n}}{x_{n+1}-1},-1\right)\in R \;\text{for}\; x\in \mathbb{S}^{n}\backslash \{p\}.\]
Denote as $\varrho$ the restirction of the standard Euclidean metric in  $\mathbb{R}^{n+1}$ to the unit sphere. $(\mathbb{S}^{n},\varrho)$ corresponds to the Alexandrow one point compactification of $R$, and $h:=s^{-1}$ is semialgebraic homeomorphism $R \longrightarrow \mathbb{S}^{n}\backslash \{p\}$. We define the metric on $\mathcal{F}_{n}$ as:
\[dist_{K}(K,L)=dist_{H}(h(K)\cup\{p\},h(L)\cup\{p\}) \quad K,L \in \mathcal{F}_{n}.\]
If just one of the sets $S,T \subset \mathbb{S}^{n}$ is empty then we put $dist_{K}(S,T)=diam\;\mathbb{S}^{n}+1$.\\
A standard reasoning shows that such a function is a metric giving the Kuratowski's convergence. Its definability follows from Proposition \ref{prop25}.
\end{proof}
\end{theorem}

To generalize the previously proven theorems for multifunctions with closed values we need to prove the following proposition:

\begin{proposition}\label{prop11}
Let $F: \mathbb{R}^{m} \longrightarrow \mathcal{P}(\mathbb{R}^{n})$ be a definable multifunction with closed values. Let $a \in dom\;F$ be a fixed point, and $\{x_{\nu}\}_{\nu=1}^{\infty}\subset \mathbb{R}^{n}$ a sequence such that $x_{\nu} \longrightarrow x_{0}$ for some $x_{0} \in \mathbb{R}^{m}$. Assume that $ dist_{K}(F(x_{\nu}),F(a)) \longrightarrow 0$ Then\\
1) there is $\limsup_{x \to x_{0}}F(x) \supset F(a).$\\
2) If $F$ is inner semi-continuous then $F(x_{0}) \subset F(a)$.
\begin{proof}
Ad 1) Let $z \in F(a)$. We keep the notation $h$ and $p$ used in the previous proof. Since \linebreak $dist_{K}(F(x_{\nu}),F(a))\longrightarrow 0$, and the distance between $h(z)$ and $p$ is constant, then $\max_{y\in h(F(x_{\nu}))}d(y,h(z)) \longrightarrow 0$. Therefore we will find a sequence $h(F(x_{\nu})) \ni z_{\nu} \longrightarrow h(z)$. Since $h$ is a homeomorphism, we get that $h^{-1}(z_{\nu}) \longrightarrow z$. Notice that $h^{-1}(z_{\nu}) \in F(x_{\nu})$, so $z \in \limsup_{x \to x_{0}}F(x)$.\\
Ad 2) Fix $y \in F(x_{0})$. From the inner semi-continuity we get that for $x_{\nu} \longrightarrow x_{0}$ there exist a sequence $F(x_{\nu})\ni y_{\nu}$ such that $y_{\nu} \longrightarrow y$. From the fact that $h$ is homeomorphism, we get that $h(y_{\nu}) \longrightarrow h(y)$. Due to the fact that $dist_{K}(F(x_{\nu}),F(a))\longrightarrow 0$, $h(y_{\nu})\in h(F(x_{\nu}))$, and $h(F(a))\cup\{p\}$ is compact, then $h(y) \in h(F(a))\cup \{p\}$. Since $h:R \longrightarrow \mathbb{S}^n\backslash \{p\}$ then $h(y)\in h(F(a))$. A homeomorphism is a bijection, so that $y \in F(a)$ which ends the proof.
\end{proof}
\end{proposition}

\begin{theorem}
The assertions of Theorems 5.13, 5.15 and 5.17 hold for definable multifunctions with closed values and with $dist_K$ replacing the Hausdorff metric.
\begin{proof}
Using Proposition \ref{prop11} we just repeat the previous proofs
\end{proof}
\end{theorem}

\section{Applications}

In the next part we will be applying the previously obtained results to certain definable multifunctions. Definitions and basic facts come from \cite{medial}. In this section we will present also the proposition concerning the closedness of the upper preimage of the multifunction $m$ of closest points.\\
Let $X\subset \mathbb{R}^{n}$ be a closed set and let $x \in \mathbb{R}^{n}$. We define the set of the closest points to $x$:
\[m(x)=\{y\in X \;|\; ||y-x||=d(x,X)\}.\]
Due to the fact that $X$ is closed, $m(x)$ is nonempty for every $x \in \mathbb{R}^{n}$. Morover there holds
\[m(x)=X \cap \mathbb{S}(x,d(x,X)),\]
so that $m(x)$ is compact as the intersection of a closed and a compact set. \\
This leads us to the definition of the medial axis which is the set of points with no unambiguity realisation of distance from $X$. This is a main concept of pattern recognition theory.

\begin{proposition}[\cite{medial} Proposition 2.17]
Let $X \subset \mathbb{R}^{n}$ be a closed definable set. Then the multifunction
\[m: \mathbb{R}^{n} \ni x \longrightarrow m(x) \in \mathcal{P}(\mathbb{R}^{n})\]
is definable and upper semi-continuous.
\begin{proof}[Proof]
Notice that:
\[\Gamma_{m}\ni(x,y) \Longleftrightarrow (y \in X) \wedge (\forall z \in X:\; ||x-y|| \leq ||x-z||) \]
which means that $\Gamma_{m}$ is a definable set, which ends the proof of definability of our function.\\
Let $y \in \limsup_{x \to a} m(x)$. From the definition $\exists \; x_{n} \rightarrow a,\; \exists\; m(x_{n})\ni y_{n} \rightarrow y$. Obviously, $y\in X$, this set being closed.
We have the estimation
\[||y-a||\leq ||a-x_{n}||+||x_{n}-y|| \leq ||a-x_{n}||+||x_{n}-y_{n}||+||y_{n}-y||=\]
\[=||a-x_{n}||+d(x_{n},X)+||y_{n}-y||\leq ||a-x_{n}||+||a-x_{n}||+d(a,X)+||y_{n}-y||.\]
For a fixed $\epsilon>0$ we will find $N \in \mathbb{N}$ such that for $n>N$, $max\{||a-x_{n}||,||y_{n}-y||\}\leq \frac{\epsilon}{3}$. Due to this,
\[||y-a|| \leq \epsilon +d(a,X)\]
so finally $y \in m(a)$.
Therefore, we get the inclusion
\[\limsup_{x \to a}m(x) \subset m(a).\]
Take $y \in m(a)$. If $a \in X$ then $a=y$ and taking any sequence $\{x_{n}\}_{n=1}^{\infty}\subset X$ such that $x_{n} \longrightarrow a$ we have $x_{n} \longrightarrow y$ because $m(x_{n})=\{x_{n}\}$. If $a \notin X$, join the points $a$ and $y$ by a segment and let $\{x_{n}\}$ be a sequence of points on this segment such that $x_{n} \longrightarrow a$. Notice that in this situation $y \in m(x_{n})$ by the strict convexity of the norm.. If there was $y \notin m(x_{n})$ then for $z \in m(x_{n})$, we would get $||y-a||=d(a,X)\leq ||a-x_{n}||+d(x_{n},X) < ||a-x_{n}||+||x_{n}-y||=d(a,X)$, what yields a contradiction. Due to this we have the reverse inclusion and finally the equality
\[\limsup_{x \to a}m(x)=m(a),\]
which ends the proof of the proposition.
\end{proof}
\end{proposition}
The above multifunction $m$ is called the multifunction of closest points. The values $m(x)$ are the closest points to $x$ in $X$.
\begin{proposition}
Let $X \subset \mathbb{R}^{m}$ be a closed set, $a \in \mathbb{R}^{n}$ a fixed point. Then the set $m_{*}(m(a))$  is closed. 
\begin{proof}[Proof] Take a sequence
$\{x_{n}\}_{n=1}^{\infty}\subset \{x \in \mathbb{R}^{m}\;|\; m(a)\subset m(x)\} $ such that $x_{n} \longrightarrow x_{0}$ for some $x_{0}\in \mathbb{R}^{n}$. We want to show that $m(a)\subset m(x_{0})$. Fix $y \in m(a)$. We have the estimation
\[||y-x_{0}|| \leq ||y-x_{n}|| +||x_{n}-x_{0}||= d(x_{n},X)+||x_{n}-x_{0}||\leq \]
\[\leq d(x_{0},X)+2||x_{n}-x_{0}||\]
For a fixed $\epsilon>0$ we will find $N>0$ such that for $n>N$, $||x_{n}-x_{0}|| \leq \frac{\epsilon}{2}$. We get the inequality $||y-x_{0}||\leq d(x_{0},X)+\epsilon$, which gives us $y \in m(x_{0})$.
\end{proof}
\end{proposition}

\begin{theorem}
Let $X \subset \mathbb{R}^{m}$ be a closed set and $a \in \mathbb{R}^{m}$. Then the set $m^{\#}(m(a))$ is closed.
\begin{proof}
Take a sequence
$\{x_{n}\}_{n=1}^{\infty}\subset m^{\#}(m(a)) $ such that $x_{n} \longrightarrow x_{0}$ for some $x_{0}\in \mathbb{R}^{n}$. This means that $m(a) \cap m(x_{n}) \neq \emptyset$. Define the following sequence:
\[ z_{n}= \text{min}_{l} \{m(x_{n})\cap m(a)\}, \;\text{ with respect to the lexicographic order.}\]
Since $m(a)$ is compact there exist $\{w_{n}\}_{n=0}^{\infty}$, a convergent subsequence of $\{z_{n}\}_{n=0}^{\infty}$ and a corresponding to it $\{x_{n}^{*}\}_{n=0}^{\infty}$, a  subsequence of $\{x_{n}\}_{n=0}^{\infty}$ such that $w_{n} \in m(x_{n}^{*})$. Denote $w$ by $\lim_{n \to \infty}w_{n}$. Obviously, $w \in m(a)$. We will show that $w \in m(x_{0})$. Fix $\epsilon>0$. Take $N>0$ such that for $n>N$: $||w_{n}-w|| \leq \frac{\epsilon}{3}$, $||x_{n}^{*}-x_{0}|| \leq \frac{\epsilon}{3}$. For $n>N$ we have the following inequalities:
 \[||w-x_{0}|| \leq ||w-w_{n}||+||x_{n}^{*}-w_{n}||+||x_{n}^{*}-x_{0}|| \leq\]
\[\leq ||w-w_{n}||+d(x_{n}^{*},X)+||x_{n}^{*}-x_{0}|| \leq\]
\[\leq ||w-w_{n}||+d(x_{0},X)+||x_{n}^{*}-x_{0}||+||x_{n}^{*}-x_{0}|| \leq\]
\[d(x_{0},X) + \epsilon.\]
Since $\epsilon$ was arbitrary it means that $w \in m(x_{0})$ and $m(x_{0}) \cap m(a) \neq \emptyset$.
\end{proof}
\end{theorem}

\begin{example}Notice that the sets $m^{*}(m(a))$ and $m^{-1}(m(a))$ may not be closed. Consider the following example. Let $X=\{1,0\} \subset \mathbb{R}$ and $a=\frac{1}{4}$. For a sequence $x_{n}=\frac{1}{2}-\frac{1}{n}$ we have that $x_{n} \longrightarrow \frac{1}{2}$ and $m(x_{n})=m(a)$ but $m(\frac{1}{2}) \not\subset m(a)$.
\end{example}

For a fixed closed semialgebraic set $X \subset \mathbb{R}^{n}$ and $a\in \mathbb{R}^{n}$ basing on the Theorem \ref{thm2} there holds
\[dist_{H}(m(x),m(a))\geq C d(x,\{y\in \mathbb{R}^{n}\;|\;m(a)\subset m(y)\})^{\alpha}, \; x \in K\]
for some $C,\alpha>0$ and $K$ being compact neighbourhood of $a$.\\
This means that for some constant $C$, the distance between the closest points to $x$ and $a$ is estimated from below by the power of the distance of $x$ to the points for which $m(a)$ are are among the closest points to $x$. \\

In a similar way we can, for a fixed closed $X$, define the set of points for which a given point is the closest point of $X$:
\[N(a)=\{x \in \mathbb{R}^{n} \;|\; a \in m(x)\}.\]
Notice that this is a non-empty set, because $a \in X$. Let $(x_{n})_{n=1}^{\infty} \subset N(a)$ be a sequence such that $x_{n} \longrightarrow x_{0}$ for some $x_{0} \in \mathbb{R}^{n}$. Then:
\[||x_{0}-a|| \leq ||x_{0}-x_{n}||+||x_{n}-a||=||x_{0}-x_{n}||+d(x_{n},X)\leq\] \[\leq ||x_{0}-x_{n}||+||x_{0}-x_{n}||+d(x_{0},X).\]
For a fixed $\epsilon >0$ we will find $N>0$ such that for $n>N$ $||x_{0}-x_{n}||\leq \frac{\epsilon}{2}$, and therefore $||a-x_{0}||\leq \epsilon +d(x_{0},X)$ what means that $x_{0} \in N(a)$, i.e. $N(a)$ is closed.

\begin{proposition}[\cite{medial} Proposition 2.2]
Let $X$ be a closed set, $a\in X$ a fixed point. Then the set $N_{*}(N(a))$ is closed.
\begin{proof}[Proof]
Let $\{x_{n}\}_{n=1}^{\infty} \subset \{x \in X\;|\; N(a)\subset N(x)\}$ be a sequence such that $x_{n} \longrightarrow x_{0}$ for some $x_{0}$. Due to the fact that $\{x_{n}\}_{n=1}^{\infty}\subset X$ we have $x_{0} \in X$. Therefore $N(x_{0})$ is well-defined. We want to show that $N(a) \subset N(x_{0})$. Fix $z \in N(a)$. 
\[||z-x_{0}||\leq ||z-x_{n}||+||x_{n}-x_{0}|| =d(z,X)+||x_{n}-x_{0}||\]
For a fixed $\epsilon>0$ we will find $N>0$ such that for $n>N$ $||x_{n}-x_{0}|| \leq \epsilon$. We get that $||z-x_{0}|| \leq d(z,X)+\epsilon$. what means that $x_{0} \in m(z)$. Then $z \in N(x_{0})$, and this give us the assertion.
\end{proof}
\end{proposition}

\begin{example}Now notice that the weak preimage of $N$ not need to be closed in general. Consider the following set $X=\{(0,0)\} \cup \{(x,y) \in \mathbb{R}^{2}\; |\; x\in [\frac{1}{2},1],\; y=-(x-1)^{2} \}$. Take a sequence $\{(x_{n},y_{n})\}_{n=1}^{\infty} \subset X \backslash \{(1,0)\}$ such that $(x_{n},y_{n}) \longrightarrow (1,0)$. Observe that $N((x_{n},y_{n})) \cap N((0,0)) \neq \emptyset$ for every $n \geq 1$, but $N((0,0)) \cap N((1,0)) = \emptyset$.
\end{example}

\begin{remark}
The weak preimage of $N$ is always closed in the one dimensional case. Let $X \subset \mathbb{R}$ be a closed set. Take $a \in X$ and a sequence $\{x_{n}\}_{n=1}^{\infty} \subset X$ tending to some $x_{0} \in \mathbb{R}$ such that $N(a) \cap N(x_{n})\neq \emptyset$. There is a subsequence $\{x_{n_{k}}\}_{k=1}^{\infty}$ with $x_{n_{k}} \leq a$ for all $k\geq1$ or there is a a subsequence $\{x_{n_{k}}\}_{k=1}^{\infty}$ with $x_{n_{k}} \geq a$ for $k\geq1$. Consider the first case and suppose there exist $k_{0},k_{1}$ such that $x_{n_{k_{0}}}<x_{n_{k_{1}}}<a$. Let $t\in \mathbb{R}$ be such that $t \in N(a)$ and $|x_{n_{k_{0}}}-t|=|a-t|$. Notice that $|x_{n_{k_{1}}}-t| <|a-t|$ and it is a contradiction with the fact that $t \in N(a)$. This means that such a sequence must be constant. We argue in the same way in the case of the second subsequence. This means that the sequence $\{x_{n}\}_{n=1}^{\infty}$ is constant so $x_{0} \in X$.
\end{remark}

\begin{proposition}[\cite{medial} Proposition 2.2]
Let $X \subset \mathbb{R}^{n}$ be a closed definable set. Then the multifunction
\[N: X \ni x \longrightarrow N(x) \in \mathcal{P}(\mathbb{R}^{n})\]
is definable and outer semi-continuous at points that are not isolated.
\begin{proof}[Proof]
Notice that:
\[\Gamma_{N}\ni(x,y) \Longleftrightarrow x \in X \wedge  \forall z \in X :\;||x-y||\leq||y-z|| \]
what means that the multifunction $N$ is definable.
Fix $a \in X$ such that $a$ is not an isolated point in $X$. Let $y \in \limsup_{x \to a}N(x)$. By definition $\exists \; x_{n} \rightarrow~a$, $\exists\; N(x_{n})\ni y_{n} \rightarrow y$. We want to show that $a$ is a closest point to $y$. Notice that we have the estimation
\[||a-y||\leq ||a-x_{n}||+||x_{n}-y_{n}||+||y_{n}-y||=||a-x_{n}||+d(y_{n},X)+||y_{n}-y||\leq\]
\[\leq ||a-x_{n}||+||y_{n}-y||+d(y,X)+||y_{n}-y||\]
For a fixed $\epsilon >0$ we will find $N>0$, such that $n>N$ $\max\{||a-x_{n}||,||y_{n}-y||\}\leq \frac{\epsilon}{3}.$ We get the estimation
\[||a-y||\leq d(y,X)+\epsilon,\]
what means that $y\in N(a)$ and we get the outer semi-continuity.
\end{proof}
\end{proposition}

For a fixed closed semialgebraic set $X \subset \mathbb{R}^{n}$ and $a\in \mathbb{R}^{n}$ that is not isolated, Proposition \ref{prop11} gives us the inequality for some $C,\alpha>0$
\[dist_{K}(N(x),N(a))\geq C d(x,\{y\in X \;|\;N(a)\subset N(y)\})^{\alpha}, x \in K\]
where $K$ is a compact neighbourhood $a$ in X.

We get that for some constant $C$, the distance between two sets for which respectively $x$ and $a$ are the closest points, is estimated from below by the power of the distance of $x$ to the set of those points of $X$ that are closest to the same points that $a$ is the closest to.\\

\end{document}